\documentclass[sn-mathphys-num]{sn-jnl}

\usepackage{graphicx}%
\usepackage{multirow}%
\usepackage{amsmath,amssymb,amsfonts}%
\usepackage{amsthm}%
\usepackage{mathrsfs}%
\usepackage[title]{appendix}%
\usepackage{xcolor}%
\usepackage{textcomp}%
\usepackage{manyfoot}%
\usepackage{booktabs}%
\usepackage{algorithm}%
\usepackage{algorithmicx}%
\usepackage{algpseudocode}%
\usepackage{listings}%

\theoremstyle{thmstyleone}%
\newtheorem{theorem}{Theorem}[section]
\theoremstyle{thmstyletwo}%

\theoremstyle{thmstylethree}%

\newtheorem{lemma}{Lemma}[section]

\begin{document}





\title[Rational methods for abstract semilinear problems without order reduction]{Rational methods for abstract semilinear problems without order reduction}

\author[1]{\fnm{Carlos} \sur{Arranz-Simón}}\email{carlos.arranz@uva.es}

\author[1]{\fnm{Bego\~na} \sur{Cano}}\email{bcano@uva.es}

\author[1]{\fnm{C\'esar} \sur{Palencia}}\email{cesar.palencia@uva.es}

\affil[1]{\orgdiv{Applied Mathematics Department}, \orgname{IMUVA,University of Valladolid}, \orgaddress{\street{P/ Belen, 7}, \city{Valladolid}, \postcode{47011},  \country{Spain}}}





\abstract{Rational methods are intended to time integrate linear homogeneous problems. However, their scope can be extended so as to cover linear nonhomogeneous problems. In this paper the integration of semilinear problems is considered. The resulting procedure requires the same computational cost than the one of a linked Runge--Kutta method, with the advantage that the order reduction phenomenon is avoided. Some numerical illustrations are included showing the predicted behaviour of the proposed methods.}
\keywords{order reduction; rational methods; Runge--Kutta methods; partial differential equations; abstract evolution equations.}


\maketitle

\section{Introduction}

It is well-known that many phenomena in science, engineering or social sciences can be described by semilinear initial value problems, which can be written in abstract form like

\begin{align}
&\left\{
\begin{aligned}
u'(t) & = A u(t) + f(t,u(t)), \quad t > 0, \\
u(0)  & = u_0,
\end{aligned}
\right.
\label{semilinearp}
\end{align}
where $A$ is a space differential operator and $f(t,u)$ is a nonlinear term which may contain derivatives in space of at least one order less than those of $A$.

On the other hand, it is also well-known that integrating these problems in time with methods with stages, as the widely used Runge--Kutta methods, lead to order reduction, i.e. the order of convergence which is observed is less than that corresponding to the same method when integrating a non-stiff ODE \cite{SSVH}. Because of this, several techniques have been devised in the literature to avoid it.

The first one consists of designing methods satisfying certain stiff (or weak stage) order conditions and not only the non-stiff ones \cite{BKRSS, BKSS, BKSS2,LV, LO}. The disadvantage of this is that there is less freedom in the choice of coefficients, which may influence the computational efficiency of the suggested methods. Moreover, in the case of weak stage order conditions, the analysis has just been performed for linear problems and just seen numerically for nonlinear problems till order $3$.

Another technique for linear problems was suggested in \cite{CP}, which consists of converting the problem, through the solution of several elliptic problems, to one for which order reduction is not observed. This procedure has the advantage to be valid for any method, but the solution of the corresponding elliptic problems also means a non-negligible computational cost. Moreover, the generalization of this technique to nonlinear problems has just been performed in \cite{CFN}, where the non-natural hypothesis that $f(t,u)$ vanishes for nul $u$ must be made.

A third procedure consists of correcting the boundary values for the stages \cite{AGC, A2, AC, AC2, P}. In the case of linear problems, if an analytic expression for the source term is known for which space and time derivatives can be calculated, order reduction can be completely avoided \cite{A2, AC}. In the case of nonlinear problems with time-dependent nonlinearity, numerical differentiation is required \cite{P} if order $\ge 3$ is pursued with Runge--Kutta methods. In the case of Rosenbrock methods, which are Runge--Kutta type integrators which make use of the Jacobian of the vector field to integrate nonlinear problems more efficiently, numerical differentiation is just required to get order $\ge 5$  if $f$ does not contain derivatives in space and to get order $\ge 4$ if it does \cite{AC2}. This technique is very cheap since just some additional calculations on the boundaries of the stages must be performed. The only disadvantage is that, in order to get high accuracy, numerical differentiation is required, which is well-known to be unstable for small grids \cite{SS}.

This paper deals with the problem of avoiding order reduction when integrating the semilinear problem (\ref{semilinearp}) by generalizing the technique already described in \cite{ArP} for linear problems. The main idea is to consider an $A(\theta)$-acceptable rational method, which can be that associated to a Runge--Kutta method when integrating a homogeneous problem, and to write suitably problem (\ref{semilinearp}) as a homogeneous one. This technique has the advantage that no numerical differentiation is required if $f$ has no space derivatives and it is also very cheap since the same number of linear systems as that corresponding to the Runge--Kutta method when integrating a linear system must be solved. On the other hand, once some starting values are calculated, just one nonlinear system per step must be calculated if the method is implicit in the part corresponding to the source term. Moreover, there is no barrier on the order of accuracy which can be achieved. The final technique resembles that suggested in \cite{K} for Garlerkin/Runge--Kutta discretizations for semilinear parabolic problems, but it is much more general in its deduction and application.

In contrast to \cite{ArP}, due to the nonlinearity, a more complex analysis must be performed distinguishing between different types of nonlinearities and considering sharp regularization estimates. On the other hand, it is not an aim of the paper to consider initial boundary value problems with time-dependent boundary conditions. For linear problems, that has already been studied in \cite{ACP} and we numerically know that it also works properly for semilinear problems, but its thorough analysis deserves further research.

The paper is structured as follows: In Section 2, the derivation of the methods is described. Then, in Section 3, some technical lemmas are given, which allow to prove convergence of the method in Section 4. As the method requires in principle to calculate $f$ at some previous values (as with multistep methods), the procedure to calculate the starting values is described and justified in Section 5. Finally, in Section 6, some numerical experiments are shown which corroborate in both hyperbolic and parabolic problems the avoidance of order reduction.

\section{Derivation of the methods}
The purpose of this paper is to extend the family of rational methods for linear problems of the form
\begin{equation}
\left\{ \begin{array}{lcl} u'(t) & = & A u(t) + f(t), \quad t > 0, \\
u(0)  & = & u_0,
\end{array} \right.
\label{linearp}
\end{equation}
which were introduced in \cite{ArP}, to semilinear problems of the form (\ref{semilinearp}). In this section we aim to explain how to do this. For the convenience of the reader, we briefly review what these methods consist of and the main assumptions we consider.

Our analysis will be based on the formulation of (\ref{semilinearp}) as an evolution equation in a Banach space $\left(X, \|\cdot \|\right)$. We assume standard properties for the operator $A$, as in \cite{Lunardi, Pazy}.

Throughout the article, we want to distinguish between the hyperbolic and parabolic case, since for the latter some results can be improved. For reasons that will become clear in the next paragraph, we will refer to the weaker hypothesis of the first case with $\alpha = 0$ and with $\alpha > 0$ to the stronger hypothesis of the second case.

\paragraph*{\textbf{Hypothesis 1.}}
\begin{list}{$\bullet$}{}
\item If $\alpha = 0$, let $A : D\left(A\right) \subset X \rightarrow X$ be a densely defined and closed linear operator on $X$ satisfying the resolvent condition
\begin{equation}
\| \left(\lambda I - A\right)^{-n}\| \leq \dfrac{M}{\left( \text{Re} \,\lambda - \omega \right)^n},
\label{resolventbd1}
\end{equation}
for $n = 1, 2, \dots$, on the plane $\left\lbrace \lambda \in \mathbb{C} :  \text{Re}\, \lambda > \omega \right\rbrace$ for $M \geq 1$, $\omega \in \mathbb{R}$.

Under this assumption, the operator $A$ is the infinitesimal generator of a $\mathcal{C}_0$ semigroup $\left\lbrace \text{e}^{tA} \right\rbrace_{t\geq 0}$ that satisfies the growth estimate
\begin{equation}
\|\text{e}^{tA}\| < M \text{e}^{t \omega}.
\label{growthsemigroup}
\end{equation}

\item If $\alpha > 0$, let $A : D\left(A\right) \subset X \rightarrow X$ be a densely defined and closed linear operator on $X$ satisfying the resolvent condition
\begin{equation}
\| \left(\lambda I - A\right)^{-1}\| \leq \dfrac{M}{\left| \lambda - \omega \right|}
\end{equation}
on the sector $\left\lbrace \lambda \in \mathbb{C} :  0 \leq \left|\text{arg}\left(\lambda - \omega\right)\right|\leq \pi - \theta, \lambda \neq \omega \right\rbrace$ for $M \geq 1$, $\omega \in \mathbb{R}$ and sectorial angle $0 < \theta < \pi/2$.

Under this assumption, the operator $A$ is the infinitesimal generator of an analytic semigroup $\left\lbrace \text{e}^{tA} \right\rbrace_{t\geq 0}$. Fixed $\omega^* \geq \omega$, the fractional powers of $\tilde{A} = \omega^* I - A$ are well defined. We set $X_{\alpha} = D(\tilde{A}^{\alpha})$ endowed with the graph norm $\|\cdot\|_{\alpha}$ of $\tilde{A}$. It is well known that $X_{\alpha}$ is independent of $\omega^* > \omega$ and that changing $\omega^* > \omega$ results in an equivalent norm. In addition to (\ref{growthsemigroup}), we now also have the estimate
\begin{equation}
\|t^{\alpha} \tilde{A}^{\alpha} \text{e}^{tA}\| \leq M \text{e}^{t \omega}.
\label{resolventbd2}
\end{equation}
\end{list}

The class of nonlinearities $f$ allowed in this setting depends on the nature of the semigroup $\left\lbrace \text{e}^{tA} \right\rbrace_{t\geq 0}$.
\paragraph*{\textbf{Hypothesis 2.}}
\begin{list}{$\bullet$}{}
\item If $\alpha = 0$ and $\left\lbrace \text{e}^{tA} \right\rbrace_{t\geq 0}$ is just a $\mathcal{C}_0$ semigroup, we assume that $f : [0,T] \times X \rightarrow X$ is globally Lipschitz continuous. Thus, there exists a real number $L$ such that
    \begin{equation}
    \label{Lipschitz1}
    \|f(t, \xi) - f(t,\eta) \| \leq L \|\xi - \eta\|
    \end{equation}
    for all $t \in [0,T]$.

    \item If $\alpha > 0$ and $\left\lbrace \text{e}^{tA} \right\rbrace_{t\geq 0}$ is analytic, we can afford stronger nonlinearities and we assume that $f : [0,T] \times X_{\alpha} \rightarrow X$ is globally Lipschitz. Thus, there exists a real number $L$ such that
    \begin{equation}
    \label{Lipschitz2}
    \|f(t, \xi) - f(t,\eta) \| \leq L \|\xi - \eta\|_{\alpha}
    \end{equation}
    for all $t \in [0,T]$.
\end{list}

We note that for the convergence proofs in the chapter, it is sufficient that (\ref{Lipschitz1}) and (\ref{Lipschitz2}) hold in a strip along the exact solution. Although, for simplicity, we assume that f is globally Lipschitz.

The methods designed in \cite{ArP} start from a rational mapping
\begin{equation}
r(z) = r_{\infty} + \sum_{\ell = 1}^k \sum_{j=1}^{m_{\ell}} r_{\ell,j} \left(1 -  w_{\ell} z\right)^{-j},
\end{equation}
that may be the stability function of a Runge--Kutta method of order $p$. Again, when the semigroup is analytic, we can consider a wider class of rational mappings as a starting point.

Let $r(z)$ be a rational mapping that approximates the exponential $\text{e}^z$ with order $p \geq 1$, that is,
\begin{equation}
r(z) - \text{e}^z = O\left(z^{p+1}\right), \qquad r(z) - \text{e}^z \neq O\left(z^{p+2}\right), \quad \text{as} \, z \rightarrow 0.
\end{equation}

\paragraph*{\textbf{Hypothesis 3.}}
\begin{list}{$\bullet$}{}
\item
If $\alpha = 0$, we assume that $r$ is A-acceptable, i.e., that $\left|r(z)\right| \leq 1$ when $\text{Re} z \leq 0$.
\item
If $ \alpha > 0$, we assume that $\gamma = \left|r_{\infty}\right| < 1$ and, if the semigroup has sectorial angle $\theta$, we also assume that $r$ is A($\vartheta$)-acceptable with $\vartheta > \theta$. This means that $\left|r(z)\right| \leq 1$ when $z$ is in $\left\lbrace z \in \mathbb{C} : \left|\text{arg}(-z)\right| < \vartheta\right\rbrace$.
\end{list}
Under these assumptions, there exists a number $\tau_0 > 0$ such that the operator
\begin{equation}
r(\tau A)  = r_{\infty}\,I + \sum_{\ell = 1}^k \sum_{j=1}^{m_l} r_{\ell,j} \left(I - \tau w_{\ell} A\right)^{-j}
\end{equation}
is well defined for every $ 0 < \tau < \tau_0$.

Basic convergence and stability results for our analysis in relation to these operators can be found in \cite{BrTh}. In this paper, the following stability bound is stated
\begin{equation}
\label{stability}
\|r^n\left(\tau A\right)\| \leq M \, C_s(n)\, \text{e}^{\omega^+ t_n}, \quad t_n = \tau \, n, \quad C_s = O\left(n^{\sigma} \right),
\end{equation}
where $0 \leq \sigma \leq 1/2$ with $\omega^+ = \max\left\lbrace 0, \omega \right\rbrace$. The stability bound is optimal (it is sharp for $A = d/dx$ in the maximum norm \cite{BrTh0}) and may be improved depending on the nature of $r$ and $A$. For our purposes, it is interesting to note that $C_s(n)$ is such that $\sigma = 0$, when one of the following is satisfied:
\begin{enumerate}
\item[(a)] $r(z) = 1/(1-z)$, which corresponds to the implicit Euler method,
\item[(b)] for $r_{m,n}(z) = P_m/Q_n$ the Padé approximant of $\text{e}^z$ with $\text{deg} (P_m) = m$ and  $\text{deg} (Q_n) = n$, whenever $n = m - 1$. This is the case of Radau methods.
\item[(c)]$X$ is a Hilbert space and $A$ an $\omega$-dissipative operator
\item[(d)] $A$ generates an analytic semigroup in a Banach space $X$.
\end{enumerate}

In \cite{BrTh} it is also proved a convergence estimate of the homogeneous problem ($f = 0$)
\begin{equation}
\label{convergence}
\|r^n\left(\tau A\right) u_0 - \text{e}^{t_n A} u_0\| \leq C_e \, M\, \,t_n \, \tau^p \, \text{e}^{\omega^+ \kappa t_n} \|A^{p+1} \, u_0\|, \qquad n \geq 1,
\end{equation}
for $u_0 \in D\left(A^{p+1}\right)$, $\kappa = \kappa(r) \geq 1$ and $C_e = C_e(r) > 0$. Note that the previous bound shows the convergence of the rational method to the solution of the homogeneous problem with order $p$, even in cases where $C_s(n)$ is not $O(1)$.

The fact that rational methods preserve the order of convergence $p$ when applied to homogeneous problems is the key starting point of our work. This is because this property is not preserved when we consider non-homogeneous problems, where the corresponding Runge--Kutta method exhibits a reduced order of convergence that has to do with the stage order of the method, rather than to $p$ itself. This phenomenon is what is known as order reduction. The optimal orders that can be achieved have been studied in \cite{IC, Ostermann, Ostermann2}. In \cite{ArP}, a brief outline of the causes of this reduction is presented.

The strategy proposed in \cite{ArP} to design the methods is the introduction of an evolution problem in a suitable product space $Z = X \times Y$, so that the first component of this system has the original problem as its solution. This is achieved by considering the space of uniformly continuous, bounded functions $Y = \mathcal{C}_{ub} \left([0,\infty),X\right)$, the operator $L : v \in Y \mapsto v(0) \in X$ and the operator $B : v \in D\left(B\right) \subset Y \mapsto v' \in Y$, whose domain is $D\left(B\right) = \left\lbrace v \in Y : v' \in Y\right\rbrace$. The operator $B$ is the infinitesimal generator of the translation semigroup, that is,
\begin{equation*}
e^{t B}v(s) = v(t+s) \quad \text{for} \, v \in Y.
\end{equation*}
Thus, the evolution problem
\begin{align}
&\left\{
\begin{aligned}
u'(t) & = A u(t)  + L v(t), \quad t \ge 0, \\
v'(t) & = \qquad \quad \ \, B v(t), \quad t \ge 0, \\
u(0)  & = u_0,\\
v(0) &  = f,
\end{aligned}
\right.
\label{Zproblem}
\end{align}
is introduced. A direct calculation proves that the first component of the solution of this problem is the solution of (\ref{linearp}). Once this is achieved, the idea is to apply the rational method to the homogeneous problem (\ref{Zproblem}), which gives as a result the approximation
\begin{align}
\begin{aligned}
\bar{u}_{n+1} & = r_{\infty} \bar{u}_{n} \\
& \quad \ + \sum_{\ell=1}^k  \sum_{j=1}^{m_l} r_{\ell,j} \left(I - \tau  w_{\ell} A\right)^{-j}\left(\bar{u}_n + \tau w_{\ell} \sum_{i=1}^j \left(I - \tau w_{\ell} A\right)^{i-1} \, L \, \left(I - \tau  w_{\ell} B\right)^{-i} r^n (\tau B) f\right) \\
& = r\left(\tau A\right) \bar{u}_n + \tau E\left(\tau A\right)\, r^n (\tau B) f, \label{abstractscheme}
\end{aligned}
\end{align}
for a certain operator $E\left(\tau A\right) : Y \rightarrow X$. This  scheme has the disadvantage that the resolvents of $B$ are more difficult to compute than the action of the semigroup itself. This is solved by the following approximation result (Lemma 4.3 in \cite{ArP})
\begin{equation}
\|L\left(I - \tau w_{\ell} B\right)^{-i}f(t_n + \cdot) - \boldsymbol{\gamma}_{\ell,i}^n \cdot f(t_n + \tau \boldsymbol{c}_n) \| \leq C \tau^p \|f^{(p)}(t_n + \cdot)\|_{\infty},
\label{aproxB}
\end{equation}
for certain values $\boldsymbol{\gamma}_{\ell,i}^n \in \mathbb{R}^p$ whenever $\boldsymbol{c}_n \in \mathbb{R}^p$ with $c_i \neq c_j$ for $i \neq j$. This result allows to approximate one of these resolvents by a linear combination of $p$ values of the source term. An advantage of these methods is that they leave freedom to choose the abscissae on which the source term $f$ is evaluated. We can choose different nodes $\boldsymbol{c}_n \in \mathbb{R}^p$ in every step so that the scheme uses the values
\begin{equation*}
f(t_n + \tau \boldsymbol{c}_n) = [f(t_n + \tau c_n^1), \dots, f(t_n + \tau c_n^p)]^T \in X^p
\end{equation*}
to compute the approximation of $u$ at time $t_{n+1}$. After choosing the nodes, the appropriate coefficients $\boldsymbol{\gamma}_{\ell,j}^n \in \mathbb{R}^p$ must be computed solving the $p \times p$ linear Vandermonde systems
\begin{equation*}
\begin{pmatrix}
1   & \cdots & 1 \\
c_n^1 &   & c_n^p \\
\vdots & \ddots  & \vdots \\
\left(c_n^1\right)^{p-1} & \cdots & \left(c_n^p\right)^{p-1}
\end{pmatrix} \boldsymbol{\gamma}_{\ell,j}^n = \begin{pmatrix}
0! F^0_{\ell, j} \\
1! F^1_{\ell, j} \\
\vdots \\
(p-1)! F^p_{\ell, j}
\end{pmatrix} ,
\end{equation*}
where the right hand side is given by the Taylor expansion
\begin{equation}
\left(1 - w_{\ell}\, z \right)^{-j} = \sum_{k=0}^{\infty} F_{\ell, j}^k \, z^k.
\end{equation}
In our implementations, only a few different $\boldsymbol{c}_n$ are used, so that they can be easily precomputed before starting the integration. After this, a step of the method is computed by
\begin{align}
u_{n+1} & = r_{\infty}\, u_n + \sum_{\ell = 1}^k \sum_{j=1}^{m_l} r_{\ell,j} \left(I - \tau w_{\ell} A\right)^{-j} \left(u_n + \tau w_{\ell}\,\sum_{i=1}^j \left(I - \tau w_{\ell} A\right)^{i-1} \, \boldsymbol{\gamma}_{\ell,j}^n \cdot f\left(t_n + \tau \boldsymbol{c}_n\right) \right) \notag \\
& = r\left(\tau A\right)u_n + \tau E_n\left(\tau A\right) f\left(t_n+\tau\boldsymbol{c}_n\right), \label{linearscheme}
\end{align}
for a certain operator $E_n(\tau A) : X^p \rightarrow X$, whose dependence on $n$ is only due to the possibility of choosing different nodes at each step. Note that the formula in the first line is written in a more explicit way (that shows how the method can be implemented), while the second line is more compact and will be useful in the analysis of the method. The main result in \cite{ArP} shows that this scheme approximates the solution of (\ref{linearp}) with order $p$, except perhaps for a reduction of 1/2 due to the stability constant $C_s(n)$. Since we are now interested in semilinear problems of the form (\ref{semilinearp}), we propose to choose integer nodes $\boldsymbol{c}_n \in \mathbb{Z}^p$, so that the times $t_n + \tau \boldsymbol{c}_n$ fall on the time grid and we can approximate the source term using the approximate values of the function by $ f(t_n,u(t_n)) \approx f(t_n,u_n)$. Although integer equispaced nodes lead to bad conditioned Vandermonde systems for high $p$ (see \cite{ACP,ArP}), that is not a problem for our purpose of avoiding order reduction (at least for $p \leq 6$, as it can be observed in the numerical experiments of Section~6).

In what follows we will consider the use of the nodes $\boldsymbol{c}_n = [-p+1,\dots,0] \in \mathbb{Z}^p$ or $\boldsymbol{c}_n = [-p+2,\dots,1] \in \mathbb{Z}^p$. The first choice requires the use of the previous values $\boldsymbol{U}_n = [u_{n-p+1}, \dots, u_n]$ to compute $u_{n+1}$, so it is explicit; whereas the second choice requires $\boldsymbol{U}_n =[u_{n-p+2}, \dots, u_{n+1}]$, and an implicit scheme turns up.

The proposed scheme can be written in a form which is analogous to (\ref{linearscheme}),
\begin{equation}
u_{n+1} = r\left(\tau A\right) u_n + \tau E_n\left(\tau A\right) f\left(t_n+\tau \boldsymbol{c}_n, \boldsymbol{U}_n\right), \quad n \geq p-1. \label{nonlinearscheme}
\end{equation}
Starting values $u_0, u_1, \cdots, u_{p-1}$, must be provided. In Section 5 we explain how to compute the first values within this framework. In Section 6 we discuss in depth the consequences of choosing each of the node possibilities.

\section{Preliminaries: discrete inequalities and regularisation}
In this section we state some results which are required to prove the convergence of the scheme (\ref{nonlinearscheme}) in the following sections.

The first lemmas are aimed at proving a variant of the discrete Gronwall lemma which is necessary for the proof of the main result of the article. The following lemma collects some bounds whose proof is elementary, but which are stated together for the sake of clarity.

\begin{lemma} Let $1 \leq k \leq n-1$, $p \geq 1$, $m^+(k) = \inf(n-1,k+p-1)$ and $\delta, \alpha \in (0,1)$. Assume that $\tau >0$ and that $t_m = m\tau$, for $0 \leq m \leq n$. Then the following inequalities hold:
\begin{align}
\sum_{m=k}^{m^+(k)} t_{n-m}^{-\alpha}  & \leq p^{1+\alpha} \, t_{n-k}^{-\alpha}, \label{bsum1} \\
\tau \sum_{m=0}^{n-1} t_{n-m}^{-\alpha} & \leq \dfrac{t_n^{1-\alpha}}{1 - \alpha}, \label{bsum2} \\
\sum_{m=k}^{m^+(k)} \delta^{n-m-1} & \leq p\,\delta^{1-p} \delta^{n-k-1}, \label{bsum3} \\
\sum_{m=1}^{n-1} \delta^{n-m-1} & \leq \dfrac{1}{1-\delta}. \label{bsum4}
\end{align}
\label{bsum}
\end{lemma}
\begin{proof}.
To prove the first inequality notice that, for $k \leq m  \leq m^+(k)$,
\begin{equation*}
\dfrac{t_{n-m}^{-\alpha}}{t_{n-k}^{-\alpha}} = \dfrac{(n-k)^{\alpha}}{(n-m)^{\alpha}} = \left(1 + \dfrac{m-k}{n-m}\right)^{\alpha} \leq (1+ p-1)^{\alpha} = p^{\alpha},
\end{equation*}
which proves (\ref{bsum1}), since the sum has at most $p$ terms. For the second inequality, notice that
\begin{equation*}
\tau \sum_{m=0}^{n-1} t_{n-m}^{-\alpha} = \tau \sum_{m=1}^n t_m^{-\alpha} \leq \tau^{1-\alpha} \int_0^n \dfrac{ds}{s^{\alpha}} \leq \dfrac{t_n^{1-\alpha}}{1-\alpha}.
\end{equation*}
The third one is true since, for $m \leq k \leq m^+(k)$,
\begin{equation*}
\delta^{n-m-1} = \delta^{k-m}\delta^{n-k-1} \leq \delta^{1-p} \delta^{n-k-1},
\end{equation*}
and again the sum has at most $p$ terms. The last inequality is just the sum of a geometric series.
\end{proof}

We now state the following lemma, which is a variant of the discrete Gronwall lemma.

\begin{lemma} Let $\tau > 0$, $N \geq 1$ and $t_n = n\tau$, $0 \leq n \leq N$. Let $\xi_n$ be a sequence of real positive numbers with $\xi_0 = 0$ and
\begin{equation}
(a)\  \xi_n^p = \sum_{k=n-p+1}^{n} \xi_k\quad  \text{or} \quad (b) \ \xi_n^p = \sum_{k=n-p+2}^{n+1} \xi_k, \quad  \text{for} \quad p-1 \leq n \leq N-1.
\label{twosupremums}
\end{equation}
Assume that there exist $\alpha\in (0,1)$, $\delta \in [0,1)$ and $K_0, K_1 \geq 0 $ such that
\begin{equation}
\max_{0\leq k\leq p-1} \xi_k \leq K_0
\label{maxGronwall}
\end{equation}
and that
\begin{equation}
\xi_{n+p-1} \leq K_0 + K_1 \sum_{m=0}^{n-1} \left(\tau\, t_{n-m}^{-\alpha} + \tau^{1-\alpha} \delta^{n-m-1}\right) \xi_{m+p-1}^p, \quad n\geq 0.
\end{equation}
Then, there exists a constant $K\geq 0 $ depending on $\gamma, \alpha, T = N\tau, K_1, p, \delta$ such that
\begin{equation}
\xi_n \leq K\,K_0 \qquad \text{for} \quad n \geq p-1.
\end{equation}
\label{GronwallLemma}
\end{lemma}

\begin{proof}.
We first assume case (a) in (\ref{twosupremums}). Notice that
\begin{align}
\xi_{n+p-1} & \leq K_0 + K_1 \sum_{m=0}^{n-1} \left(\tau \, t_{n-m}^{-\alpha} + \tau ^{1-\alpha} \delta^{n-m-1}\right) \sum_{k=m-p+1}^m \xi_{k+p-1} \notag \\
& = K_0 + K_1 \sum_{k=-p+1}^{n-1} \sum_{m = m^-(k)}^{m^+(k)} \left(\tau t_{n-m}^{-\alpha} + \tau^{1-\alpha} \delta^{n-m-1}\right) \xi_{k+p-1}, \notag
\end{align}
where $m^-(k) = \max\left\lbrace k,0 \right\rbrace$ and $m^+(k) = \min \left\lbrace k+p-1,n-1 \right\rbrace$ are the values that allow the previous sum to be reordered. Now, we use the estimates in Lemma \ref{bsum}. For $-p+1 \leq k \leq 0$, $m^-(k) = 0$, and formulae (\ref{bsum2}), (\ref{bsum4}) imply that, for a constant $K > 0$ depending on $T, \delta, \alpha,$ it is true that
\begin{equation*}
\sum_{m=0}^{m^+(k)} \left(\tau\, t_{n-m}^{-\alpha} + \tau^{1-\alpha} \, \delta^{n-m-1}\right) \xi_{k+p-1} \leq K \, \xi_{k+p-1}.
\end{equation*}
On the other hand, for $1 \leq k \leq n-1$, $m^-(k) = k$, and taking into account (\ref{bsum1}), (\ref{bsum3}),
\begin{equation*}
\sum_{m=k}^{m^+(k)} \left(\tau\, t_{n-m}^{-\alpha} + \tau^{1-\alpha} \, \delta^{n-m-1}\right) \xi_{k+p-1} \leq K \left(\tau \, t_{n-k}^{-\alpha} + \tau^{1-\alpha} \delta^{n-k-1}\right) \xi_{k+p-1}.
\end{equation*}
Then, we combine the latter and (\ref{maxGronwall}) to get that, for $n \geq 0$,
\begin{equation}
\xi_{n+p-1} \leq K \, K_0 + K K_1 \sum_{k=1}^{n-1} \left(\tau \, t_{n-k}^{-\alpha} + \tau^{1-\alpha} \delta^{n - k - 1} \right) \xi_{k+p-1},
\label{GronwallPr}
\end{equation}
for another constant $K$. If we consider case (b), we obtain an additional term $K\, K_1\, \tau^{1-\alpha} \xi_{n+p-1}$ in the right hand side. It is clear that, for small enough $\tau$, case (b) may be reduced to formula (\ref{GronwallPr}). The proof concludes applying Lemma 2.1. in \cite{Gonzalez}.
\end{proof}

Hereafter, the letter $K$ denotes general positive constants that may depend on the semigroup ($M$, $\omega$, $\alpha$), the rational method ($C$, $\gamma$) or the interval $[0,T]$ of integration, but that does not depend on any considered particular solution $u$, source term $f$ or step-size $\tau$.

When $A$ generates an analytic semigroup and we work with functions $u:[0,T] \rightarrow X_{\alpha}$, we expect the numerical approximations to the solution to be in the space $X_{\alpha}$, not just in $X$. Since the nonlinearity $f$ takes values in $X$, the linear part of the numerical scheme, governed by the operator constructed from the rational function $r\left(\tau A \right)$, must have some regularisation property that guarantees that the numerical solutions are in $X_{\alpha}$. This is what motivates the results with which this section ends.

Notice that $r\left(\tau A\right) - r_{\infty} I$ is a linear combination of powers of resolvents of $A$. This implies, by (\ref{resolventbd1}), that for $0 < \tau < \tau_0$,
\begin{equation*}
\|\left(r\left(\tau A\right) - r_{\infty} I\right)x\| \leq K \|x\|, \quad x \in X,
\end{equation*}
and in the analytic case, Lemma 2.2 in \cite{Gonzalez} also implies that
\begin{equation*}
\|A \left(r\left(\tau A\right) - r_{\infty} I\right)x\| \leq \dfrac{K}{\tau} \|x\|, \quad x \in X,
\end{equation*}
for another constant $K > 0$ that may depend on $r$ and $T$. Then, by interpolation (see e.g. \cite{Triebel}), we get that
\begin{equation}
\|\left(r\left(\tau A\right) - r_{\infty} I\right)x\|_{\alpha} \leq \dfrac{K}{\tau^{\alpha}} \|x\|, \quad x \in X.
\end{equation}
This shows that the linear part of the numerical scheme regularises the solution after one step. For several steps, the equation (7) in \cite{Gonzalez}, although stated in a more general framework of variable step size problems, can be realised in our context as a generalisation of the above to several steps,
\begin{equation}
\|\left(r^n\left(\tau A\right) - r_{\infty}^n I\right)x\|_{\alpha} \leq \dfrac{K}{t_n^{\alpha}} \|x\|, \quad x \in X.
\label{ratregularizat}
\end{equation}
Taking into account formulae (\ref{abstractscheme}) and (\ref{linearscheme}), it can be seen that the operators $E(\tau A)$ and $E_n(\tau A)$ are given by a linear combination of resolvents of $A$. Therefore, using the same argument, it can be directly proved  that, for $0 \leq \beta \leq \alpha$,
\begin{align}
\|E(\tau A)\, v\|_{\beta} & \leq K \, \tau^{-\beta} \|v\|_Y, \quad \text{for} \, v \in Y, \label{regularisationE} \\
\|E_n(\tau A)\, v_p\|_{\beta} & \leq K \, \tau^{-\beta} \|v_p\|_{X^p}, \quad \text{for} \, v_p \in X^p, \label{regularisationEn}
\end{align}
where $\|\cdot \|_{X^p}$ corresponds to the maximum of the norm of each component in $X$.

To conclude, we state a lemma which is based on these results that will be useful in the proof of the main theorem.

\begin{lemma}\label{lm-rsum}
Let $0 < \alpha < 1$. Under hypotheses 1 and 3, let $\xi_m \in X_{\alpha}$, $0 \leq m \leq n$ and $0 < \tau < \tau_0$. Then, there exists a positive constant $K$ (that may be different in each case) such that the following estimates hold
\begin{align}
\left\| \tau \sum_{m = 0}^{n-1} r^{n-m-1}\left(\tau A\right) \xi_m \right\|_{\alpha} & \leq  K \tau \sum_{m=0}^{n-2} \left(\dfrac{\|\xi_m\|}{t_{n-m-1}^{\alpha}}  +  \gamma^{n-m-1} \|\xi_m\|_{\alpha} \right) + \tau \| \xi_{n-1} \|_{\alpha},\label{r-sum} \\
\left\| \tau \sum_{m = 0}^{n-1} r^{n-m-1}\left(\tau A\right) \xi_m \right\|_{\alpha} & \leq  K \left(\max_{0 \leq m \leq n-2} \|\xi_m\| + \tau \, \max_{0 \leq m \leq n-1} \|\xi_m\|_{\alpha} \right).
\label{r-sum2}
\end{align}
\end{lemma}
\begin{proof}.
Taking into account the regularization estimate (\ref{ratregularizat}), the left hand side in (\ref{r-sum}) and (\ref{r-sum2}) is bounded by
\begin{align*}
&  \left\| \tau \sum_{m = 0}^{n-1} \left(r^{n-m-1}\left(\tau A\right) - r_{\infty}^{n-m-1}\right) \xi_m \right\|_{\alpha}  + \quad \left\| \tau \sum_{m = 0}^{n-1} r^{n-m-1}_{\infty} \xi_m \right\|_{\alpha} \\
& \leq  \tau \sum_{m=0}^{n-2} \dfrac{K}{t_{n-m-1}^{\alpha}}\|\xi_m\| + \tau \|\xi_{n-1}\|_{\alpha} + \tau \sum_{m=0}^{n-2} \gamma^{n-m-1} \|\xi_m\|_{\alpha},
\end{align*}
which proves (\ref{r-sum}). To prove (\ref{r-sum2}), notice that
\begin{align*}
\tau \sum_{m=0}^{n-2} \dfrac{\| \xi_m \|}{t_{n-m-1}^{\alpha}} & \leq \tau^{1-\alpha} \left(\int_0^n \dfrac{1}{s^{\alpha}}\, ds \right) \max_{0 \leq m \leq n-2} \| \xi_m \| \\
& \leq \dfrac{\tau^{1-\alpha}n^{1-\alpha}}{1-\alpha}\, \max_{0\leq m\leq n-2} \|\xi_m\| \leq \dfrac{T^{1-\alpha}}{1-\alpha}\, \max_{0\leq m\leq n-2} \|\xi_m\|,
\end{align*}
and
\begin{align*}
\tau \sum_{m=0}^{n-1} \gamma^{n-m-1} \|\xi_m\|_{\alpha} \leq \tau \left(\sum_{m=0}^{\infty} \gamma^m\right) \max_{0\leq m \leq n-1} \|\xi_m\|_{\alpha} \leq \tau \dfrac{1}{1-\gamma} \max_{0\leq m \leq n-1} \|\xi_m\|_{\alpha}.
\end{align*}
\end{proof}
In the case $\alpha = 0$, instead of the above lemma it will be sufficient to use the direct bound
\begin{equation}
\left\| \tau \sum_{m = 0}^{n-1} r^{n-m-1}\left(\tau A\right) \xi_m \right\| \leq \tau \,C_s(n)\, \sum_{m=0}^{n-1} \|\xi_m\|.
\label{r-sum3}
\end{equation}

\section{Convergence of the method}
Before stating the main theorem, we have to prove an auxiliar result, which is an extension of the convergence theorem (Theorem 4.4 in \cite{ArP}) to the framework of spaces $X_{\alpha}$. Its proof is similar to the one of that theorem, now taking into account the regularization estimates.

For the rest of the section, assume hypotheses 1, 2 and 3, and let $h : [0,T] \rightarrow X$ be $h(t) = f(t,u(t))$. The linear problem
\begin{align}
&\left\{
\begin{aligned}
v'(t) & = A v(t) + h(t), \quad t \ge 0, \\
v(0)  & = u_0,
\end{aligned}
\right.
\label{2linearp}
\end{align}
has $u$ as a solution and is now discretised by means of the recurrence
\begin{align}
v_{n+1} = r\left(\tau A\right)v_n + \tau E_n\left(\tau A\right) h\left(t_n+\tau\boldsymbol{c}_n\right), \quad n \geq 1, \label{2linearscheme}
\end{align}
for some sequence $\left\lbrace \boldsymbol{c}_n \right\rbrace_{n=0}^{N-1}$.

\begin{theorem}
Under hypotheses of Lemma \ref{lm-rsum}, let $u:[0,\infty) \rightarrow X_{\alpha}$ be the solution of (\ref{2linearp}) to be approximated on the interval $[0,T]$ with constant step-size $0 < \tau = T/N < \tau_0$. Assume also that $u \in \mathcal{C}^{p+1} \left([0,T], X_{\alpha} \right)$, $h \in \mathcal{C}^{p+1} \left([0,T], X\right)$. If $v_n$ is the numerical approximation to $u(t_n)$ given by (\ref{2linearscheme}),
\begin{equation}
\left\|u(t_n) - v_n\right\|_{\alpha} \leq K \, \tau^p \left(\|u^{(p+1)}\|_{\alpha,\infty} + \|h^{(p)}\|_{\infty} + \|h^{(p+1)}\|_{\infty}\right), \qquad 0 \leq n \leq N.
\label{linearconvergence}
\end{equation}
\label{th-linearconvergence}
\end{theorem}
\begin{proof}.
Using the notation of \cite{ArP}, it is straightforward to prove that $S_G$ is a semigroup in $X_{\alpha} \times Y$ and the main result in \cite{BrTh} guarantees that, for the abstract scheme
\begin{equation}
\quad \bar{v}_{n+1} = r\left(\tau A \right) \bar{v}_n + \tau E(\tau A)\, r^n\left(\tau B\right)h, \quad n \geq 1,
\label{exactscheme}
\end{equation}
we get global error of order $p$,
\begin{equation}
\|u(t_n) - \bar{v}_n \|_{\alpha} \leq C\, \tau^p \left(\|u^{(p+1)}\|_{\alpha,\infty} + \|h^{(p+1)}\|_{\infty}\right).
\end{equation}
Then subtracting (\ref{exactscheme}) from (\ref{2linearscheme}),
\begin{eqnarray*}
v_{n+1} - \bar{v}_{n+1} & = & r\left(\tau A\right) \left(v_n - \bar{v}_n\right) + \tau \left(E_n(\tau A)\,h(t_n+ \tau \boldsymbol{c}_n) - E(\tau A)\, r^{n}\left(\tau B\right)h\right) \\
& = & r\left(\tau A\right) \left(v_n - \bar{v}_n\right) + \tau \left(E_n(\tau A)\,h(t_n+ \tau \boldsymbol{c}_n) - E(\tau A)\,h(t_n + \cdot)\right)\\
&  & + \,\,\, \tau \left(E(\tau A) \left(h(t_n + \cdot) - r^{n}\left(\tau B\right)h \right)\right), \quad \text{for} \ n \geq 0,
\end{eqnarray*}
with $v_0 = \bar{v}_0$. Then, by the variation-of-constants formula, the error $\|v_n - \bar{v}_n\|_{\alpha}$ is bounded by the sum of the two terms
\begin{eqnarray*}
(I) & = & \left\|\tau \sum_{m=0}^{n-1}  r^{n-m-1}\left(\tau A\right) \left( E_m(\tau A)\, h(t_m + \tau \boldsymbol{c_m}) - E(\tau A)\, h(t_m + \cdot)\right) \right\|_{\alpha}, \\
(II) & = & \left\|\tau \sum_{m=0}^{n-1}  r^{n-m-1}\left(\tau A\right)  E(\tau A) \left(h(t_m + \cdot) - r^{m}\left(\tau B\right)h \right)\right\|_{\alpha}.
\end{eqnarray*}
The proof is concluded taking into account (\ref{r-sum2}), the regularisation  estimates (\ref{regularisationEn}) and (\ref{regularisationE}), and the approximation estimate (\ref{aproxB}). Notice that in this case $C_s(n) = O(1)$, because the semigroup is analytic.
\end{proof}

Now, we are in position to state and prove the main result.

\begin{theorem}
For $0\le \alpha <1$, let $u: [0,T] \rightarrow X_{\alpha}$ be the solution of (\ref{semilinearp}) to be approximated in the interval $[0,T]$ with constant step size $0 < \tau = T/N < \tau_0$. Let us assume hypotheses 1, 2 and 3 and also that $u \in \mathcal{C}^{p+1}\left([0,T], X_{\alpha}\right)$ and $h \in \mathcal{C}^{p+1}\left([0,T], X\right)$. If $u_n$ is the numerical approximation to $u(t_n)$ given by (\ref{nonlinearscheme}), and $u_0, \cdots, u_{p-1} \in X_{\alpha}$ are starting values satisfying
\begin{equation}
\| u(t_n) - u_n \|_{\alpha} \leq C_0 \, \tau^p, \qquad 0 \leq n \leq p-1,
\label{startingvalues}
\end{equation}
then,
\begin{equation}
\|u(t_n) - u_n \|_{\alpha} \leq K \, C_s(n) \, \tau^p \left(\|u^{(p+1)}\|_{\alpha,\infty} + \|h^{(p)}\|_{\infty} + \|h^{(p+1)}\|_{\infty} \right), \qquad 0 \leq n \leq N.
\end{equation}
\label{mainth}
\end{theorem}
\begin{proof}.
Along the proof we denote $\mathbf{f}_n = f(t_n + \tau \boldsymbol{c}_n,\mathbf{U}_n)$, $\mathbf{h}_n = h(t_n + \tau \boldsymbol{c}_n)$ and $e_n = \|u(t_n) - u_n \|_{\alpha}$, for $0 \leq n \leq N$. Set
\begin{equation}
e_n^p =
\left\{ \begin{array}{lcl} \displaystyle \sum_{k=n-p+1}^n e_k, \quad \text{if} \ \boldsymbol{c}_n = [-p+1, \cdots, 0], \\
\displaystyle \sum_{k=n-p+2}^{n+1} e_k, \quad \text{if} \ \boldsymbol{c}_n = [-p+2, \cdots, 1],
\end{array} \right.
\end{equation}
for $p-1 \leq n \leq N$. We recall (\ref{nonlinearscheme}) and (\ref{2linearscheme}) to get, for $0 \leq n \leq N$,
\begin{equation}
u_{n+p} - v_{n+p}  =  r\left(\tau A\right) \left(u_{n+p-1} - v_{n+p-1} \right) +  \tau E_{n+p-1}(\tau A) \left(\mathbf{f}_{n+p-1} - \mathbf{h}_{n+p-1}\right).
\end{equation}
By the discrete variation-of-constants formula,
\begin{align}
u_{n+p-1} - v_{n+p-1} = & r\left(\tau A \right)^n \left(u_{p-1} - v_{p-1} \right) \notag \\
& + \tau \sum_{m=0}^{n-1} r\left(\tau A \right)^{n-m-1} E_{m+p-1}(\tau A) \left(\mathbf{f}_{m+p-1} - \mathbf{h}_{m+p-1}\right), \label{vocnonl}
\end{align}
for $0 \leq n \leq N-p+1$. We use (\ref{regularisationEn}) and the Lipschitz property of $f$ to get
\begin{equation*}
\| E_{m+p-1}(\tau A) \left(\mathbf{f}_{m+p-1} - \mathbf{h}_{m+p-1} \right) \| \leq K\, L \, e_{m+p-1}^p, \qquad 0 \leq m \leq N - p,
\end{equation*}
and
\begin{equation*}
\| E_{m+p-1}(\tau A) \left(\mathbf{f}_{m+p-1} - \mathbf{h}_{m+p-1} \right) \|_{\alpha} \leq \tau^{-\alpha} \,K\, L \, e_{m+p-1}^p, \qquad 0 \leq m \leq N - p.
\end{equation*}
Then, we bound the sum in (\ref{vocnonl}) by combining the previous estimates together with (\ref{r-sum}). On the other hand, the first term in (\ref{vocnonl}) is bounded using (\ref{stability}), (\ref{linearconvergence}) and (\ref{startingvalues}), giving rise to
\begin{align*}
\|u_{n+p-1} - v_{n+p-1} \|_{\alpha} \leq M \, C \, C_s(n)\, \tau^p  \left( \|u^{(p+1)}\|_{\alpha,\infty} + \|h^{(p)}\|_{\infty} + \|h^{(p+1)}\|_{\infty}\right) \\
+ K\,L\,C_s(n)\,\sum_{m=0}^{n-1} \left( \tau \, t_{n-m}^{-\alpha} + \tau^{1-\alpha} \, \gamma^{n-m-1} \right)\, e_{m+p-1}^p,
\end{align*}
for some other constants $C$ and $K$. Notice that, due to (\ref{r-sum3}), if $\alpha = 0$ the term $\tau^{1-\alpha} \, \gamma^{n-m-1}$ is unnecessary whereas if $\alpha >0$ then $C_s(n) = O(1)$. Finally, to bound the global error we combine the above estimate with Theorem \ref{th-linearconvergence} to get
\begin{align*}
e_{n+p-1}  \leq & M\,C\,C_s(n)\, \tau^p \, \left( \|u^{(p+1)}\|_{\alpha, \infty} + \|h^{(p)}\|_{\infty} + \|h^{(p+1)}\|_{\infty}\right) \\ & + K\,L\,C_s(n)\, \sum_{m=0}^{n-1} \left( \tau \, t_{n-m}^{-\alpha} + \tau^{1-\alpha} \, \gamma^{n-m-1} \right)\, e_{m+p-1}^p,
\end{align*}
and the proof concludes by using the version of discrete Gronwall lemma in Lemma \ref{GronwallLemma}. Notice that the hypothesis (\ref{startingvalues}) is fully taken into account in this step.
\end{proof}

\section{Starting values}
The scheme which has been presented requires evaluating the source term at each time step. To do so, previously computed approximated values $u_n$ can be used to evaluate $f$ at $t = t_n$. Even so, we require some starting values $u_0, u_1, \dots, u_{p-1}$ to compute the first step and start the recurrence process.

One first possibility is just to use an auxiliary method to compute the starting values. However, in this context it is natural to look for values $u_0, \dots, u_{p-1}$ that satisfy the implicit scheme
\begin{equation}
\label{ImpScheme}
u_{n+1} = r\left(\tau A\right)u_n + \tau E_n\left(\tau A\right)\, f\left(t_n + \tau \boldsymbol{c}_n, \mathbf{U}_n \right), \quad 0 \leq n \leq p-2,
\end{equation}
where in the first steps $\boldsymbol{c}_n$ is such that $t_n + \tau \boldsymbol{c}_n = [0, \tau, \dots, (p-1)\tau]^T$ and $\mathbf{U}_n = [u_0, \dots, u_{p-1}]^T$ for $0 \leq n \leq p-2$. It is then necessary to show that the system (\ref{ImpScheme}) has a unique solution that approximates the values $u(t_n)$, $1\leq n \leq p-1$, within the adequate order. To see this, we rewrite the system as a fixed point equation in $X_{\alpha}^{p-1}$
\begin{equation}
U^* = \mathcal{N}\left(U^*\right),
\label{fixedpoint}
\end{equation}
where $\mathcal{N}: X_{\alpha}^{p-1} \rightarrow X_{\alpha}^{p-1}$ is the mapping defined by
\begin{equation*}
\mathbf{U} = \begin{pmatrix}
u_1 \\
\vdots \\
 \\
u_{p-1}
\end{pmatrix} \mapsto \mathcal{N}(\mathbf{U}) =
\begin{pmatrix}
\tilde{u}_1 \\
\vdots \\
 \\
\tilde{u}_{p-1}
\end{pmatrix},
\end{equation*}
where the $\tilde{u}_j$, $0 \leq j \leq p-1$, are defined recursively by $\tilde{u}_0 = u_0$ and
$$
\tilde{u}_{j+1} = r\left(\tau A\right) \tilde{u}_{j} + \tau E_n(\tau A) \, f\left(t_j + \tau \boldsymbol{c}_{j}, [\tilde{u}_0, \tilde{u}_1, \dots, \tilde{u}_j,u_{j+1},\dots,u_{p-1}]^T\right),
$$
for $0 < j < p-1$. To show that (\ref{fixedpoint}) has a unique solution it suffices to see that $\mathcal{N}$ is a contractive mapping for sufficiently small $\tau$. In fact, if $\mathbf{U}, \mathbf{V} \in X_{\alpha}^{p-1}$, the first component of $\mathcal{N}(\mathbf{U}) - \mathcal{N}(\mathbf{V})$ is
\begin{equation}
\tilde{u}_1 - \tilde{v}_1 = \tau E_n(\tau A) \, \left( f\left(t_0 + \tau \boldsymbol{c_0},[u_0, \mathbf{U}^T]^T\right) - f\left(t_0 + \tau \boldsymbol{c_0},[u_0, \mathbf{V}^T]^T\right)\right),
\end{equation}
so using (\ref{regularisationEn}) and the Lipschitz property of $f$
\begin{equation}
\|\tilde{u}_1 - \tilde{v}_1\|_{\alpha} \leq K L\tau^{1-\alpha} \|\mathbf{U}- \mathbf{V} \|_{\alpha},
\end{equation}
where $\|\mathbf{U}\|_{\alpha}$ for a vector $\mathbf{U} \in X_{\alpha}^p$ denotes the maximum of the $\|\cdot\|_{\alpha}$-norm of its components. Then we assume that $\|\tilde{u}_k - \tilde{v}_k\|_{\alpha} \leq k\,K^k\,L\tau^{1-\alpha} \|\mathbf{U}- \mathbf{V} \|_{\alpha}$ and proceed by induction for $1 \leq k \leq p-2$,
\begin{align*}
\|\tilde{u}_{k+1} - \tilde{v}_{k+1}\|_{\alpha} & \leq K \|\tilde{u}_k - \tilde{v}_k\|_{\alpha} + K L\tau^{1-\alpha} \max\left\lbrace \max_{1 \leq j \leq k} \|\tilde{u}_j - \tilde{v}_j \|_{\alpha}, \|\mathbf{U}- \mathbf{V} \|_{\alpha} \right\rbrace \\
& \leq(k+1) \, K^{k+1}\,L\,\tau^{1-\alpha} \|\mathbf{U}- \mathbf{V} \|_{\alpha}.
\end{align*}
Therefore,
\begin{equation}
\|\mathcal{N}(\mathbf{U}) - \mathcal{N}(\mathbf{V}) \|_{\alpha} \leq p\,K^p\, L \, \tau^{1-\alpha} \|\mathbf{U} - \mathbf{V}\|_{\alpha},
\label{contaction}
\end{equation}
and the mapping is contractive for sufficiently small $\tau$. In that case, the contractive mapping theorem guarantees that (\ref{fixedpoint}) has a unique solution $ \mathbf{U} = [u_1 \dots,u_{p-1}]^T$. To conclude, we show that this fixed point approximates the solution with order $p$. We set $\mathbf{U}(t_n) = [u_0, u(t_1), \dots, u(t_{p-1})]^T$ for $0 \leq n \leq p-2$. Under the assumptions of Theorem \ref{th-linearconvergence}, the scheme
\begin{equation}
\bar{u}_{n+1} = r\left(\tau A\right) \bar{u}_n + \tau E_n(\tau A) f\left(t_n + \tau \boldsymbol{c}_n, \mathbf{U}(t_n)\right), \quad 0 \leq n \leq p-2,
\end{equation}
is such that
\begin{equation}
\|\bar{u}_n - u(t_n) \|_{\alpha} \leq C\, K \, \tau^p\left( \|u^{(p+1)}\|_{\alpha,\infty} + \|h^{(p)}\|_{\infty} + \|h^{(p+1)}\|_{\infty}\right), \quad 0 \leq n \leq p-1.
\label{bound1SV}
\end{equation}
Then we have, for $0 \leq n \leq p-2$,
\begin{equation*}
u_{n+1} - \bar{u}_{n+1} = r\left(\tau A\right)\left(u_n - \bar{u}_n \right) + \tau E_n(\tau A) \left(f(t_n + \tau \boldsymbol{c}_n,\mathbf{U}_n) - f(t_n + \tau \boldsymbol{c}_n,\mathbf{U}(t_n))\right),
\end{equation*}
and by the discrete variation-of-constants formula and (\ref{regularisationEn}), we get that $\|u_n - \bar{u}_n \|_{\alpha}$ is bounded by
\begin{align*}
 & \tau \left\|\sum_{m=0}^{n-1} r\left(\tau A \right)^{n-m-1} E_m(\tau A) \left(f(t_m + \tau \boldsymbol{c_m},\mathbf{U}_m) - f(t_m + \tau \boldsymbol{c_m},\mathbf{U}(t_m))\right)\right\|_{\alpha} \\
&\ \ \leq p K L \tau^{1-\alpha} \sup_{1 \leq k \leq p-1} \|u(t_k) - u_k \|_{\alpha}.
\end{align*}
By the triangle inequality, the previous bound and (\ref{bound1SV})
\begin{align}
\|u_n - u(t_n)\|_{\alpha} \leq & \, K \, \tau^p\left( \|u^{(p+1)}\|_{\alpha,\infty} + \|h^{(p)}\|_{\infty} + \|h^{(p+1)}\|_{\infty}\right) \notag \\
 & + p K L\,\tau^{1-\alpha} \sup_{1\leq k \leq p-1} \|u(t_k) - u_k\|_{\alpha}.
\end{align}
We can take the supremum in the left hand side and, for small enough $\tau$, and another constant $K$,
\begin{equation}
\sup_{1\leq k \leq p-1} \|u(t_k) - u_k\|_{\alpha} \leq  K \, \tau^p\left( \|u^{(p+1)}\|_{\alpha,\infty} + \|h^{(p)}\|_{\infty} + \|h^{(p+1)}\|_{\infty}\right).
\end{equation}
In practice, to solve (\ref{fixedpoint}), we can compute the initial approximation with an auxiliary method, which in our experiments is Euler implicit method, and then iterate the function $\mathcal{N}$ to obtain initial values within the adequate order.

\section{Numerical illustrations}
In this section we show numerical results which are obtained with the numerical scheme (\ref{nonlinearscheme}) and different choices of the nodes. We deal with simple PDEs which are integrated by the method of lines. For the spatial discretization, we use finite difference methods. We take $h > 0$ as discretization parameter and we are led to systems of ODEs
\begin{align}
&\left\{
\begin{aligned}
u'_h(t) & = A_h u_h(t) + f_h(t,u_h(t)), \quad t \ge 0, \\
u_h(0)  & = u_{0,h}.
\end{aligned}
\right.
\label{discretep}
\end{align}
We adjust $f_h$ in such a way that the restriction of $u$ to the discrete mesh is the exact solution of (\ref{discretep}), so that we consider only due to the time integration. The implementation of the scheme (\ref{nonlinearscheme}) applied to the latter requires evaluating the operator $r\left(\tau A_h\right)$. In practice, this means dealing with systems of equations of the form
\begin{equation*}
\left(I - \tau \, w_{\ell} \, A_h \right)x = y,
\end{equation*}
so it is sufficient to have a routine that solves such systems, which matrices typically have a sparse structure that can be exploited.

After the spatial discretization, scheme (\ref{nonlinearscheme}) is applied to the PDE in the discretized form (\ref{discretep}). In Section~6 of \cite{ArP}, accurate details on the implementation issues of the linear version (\ref{linearscheme}) can be found. Such matters are the same in this case by adding the dependence of the source term on the function $u$ and choosing the appropriate nodes. In all examples, the scheme is implemented in three ways:
\begin{enumerate}
\item \textit{Explicit mode.} We choose $\boldsymbol{c}_n = [-p+1, \dots, 0]^T \in \mathbb{R}^p$, so $u_{n-p+1}, \dots, u_n$ are used to compute $u_{n+1}$.
\item \textit{Semiexplicit mode}. First, we take a step with the explicit scheme to have an approximation $\tilde{u}_{n+1}$ to $u(t_{n+1})$. Then, we correct this approximation by taking $\boldsymbol{c}_n = [-p+2, \dots, 1]$ and using $u_{n-p+2}, \dots, \tilde{u}_{n+1}$ to compute $u_{n+1}$.
\item \textit{Implicit mode.} We set a tolerance $TOL$ and iterate the previous process until two successive iterants $\tilde{u}_{n+1}^{[k]}$ and $\tilde{u}_{n+1}^{[k+1]}$ are such that $$\|\tilde{u}_{n+1}^{[k]} - \tilde{u}_{n+1}^{[k+1]}\|_{\alpha} \leq TOL. $$
\end{enumerate}
In all numerical experiments below, the tolerance to calculate implicitly the starting values or for the iteration in the implicit mode has been $TOL= 10^{-14}$.

\paragraph*{Example 1.}  We consider a semilinear parabolic problem in the unit interval with homogeneous Dirichlet boundary conditions,
\begin{equation}
\left\{ \begin{array}{lcll} u_t(t,x) & = & u_{xx}(t,x) + \lambda \left(\displaystyle\int_0^1 u(t,x)\,dx \right) \, u_{x} + f(t,x), & \ 0\le t \le 1, \quad 0 \leq x \leq 1, \\
u(0,x)  & = & u_0(x), & \quad 0 \leq x \leq 1, \\
u(t,0) & = & 0, & \quad 0 \leq t \leq 1, \\
u(t,1) & = & 0, & \quad 0 \leq t \leq 1.
\end{array} \right.
\label{ej1}
\end{equation}
where $f : [0,1] \times X_{\alpha} \rightarrow \mathbb{C}$, $u_0 : [0,1] \rightarrow X_{\alpha}$, $\lambda > 0$. In order to fit the problem in our framework, we take $X = L^2[0,1]$, $A = d^2/dx^2$ with $D(A) = H^2[0,1] \cap H^1_0[0,1]$ and $\alpha = 1/2$, so that $\| \cdot \|_{\alpha}=\|\cdot \|_{H^1(0,1)}$. We adjust the data $u_0$ and $f$ in such a way that $u(t,x) = x(1-x)\,\text{e}^t$, $0 \leq t, x \leq 1$, is the solution of the problem. We consider various values of the parameter $\lambda$ to test how the different implementations of the scheme behave with respect to the stiffness of the source term $f$.
\begin{figure}[htpb]
\centering
\includegraphics[height=70mm]{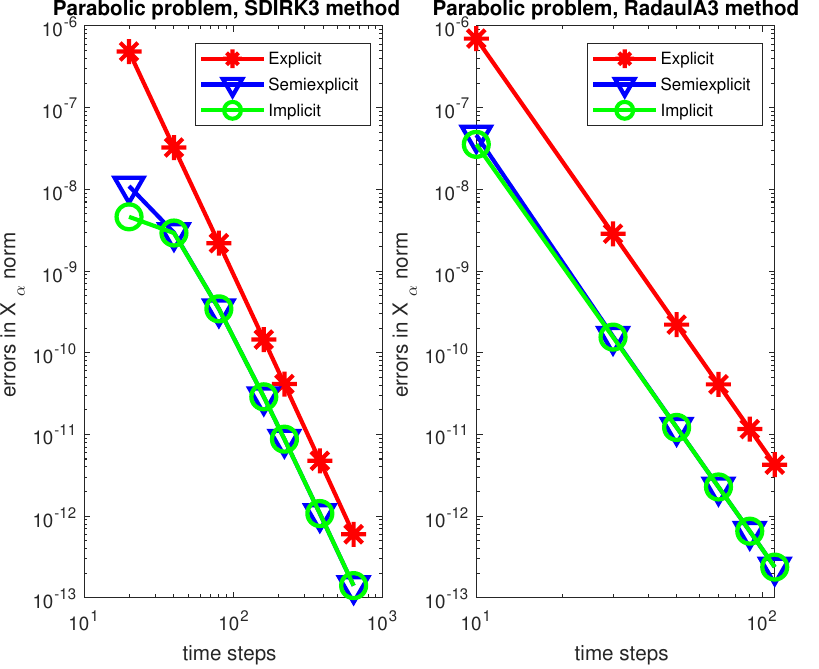}
\caption{This figure shows the error in the discrete norm $\|\cdot\|_{H^1(0,1)}$ for the parabolic problem (\ref{ej1}) with $\lambda = 1$, two different time integrators and $J = 100$.}
\label{fig1}
\end{figure}
We discretize the problem in space by means of finite differences. To this end, we fix a number $J$ of uniformly distributed nodes $x_j = jh, 1 < j < J$, in $(0, 1)$, with $h = 1/(J + 1)$. The spatial derivatives $u_x$ and $u_{xx}$ are approximated by using central finite differences and the standard three-point finite difference scheme, respectively; while the integral has been approximated by using the composite Simpson's rule. In this way, since the exact solution is a polynomial of degree two in x, there are no spatial errors.

\begin{table}[htpb]
\centering
\caption{Errors and order of convergence in Example 1 with the rational SDIRK3 method.}
\begin{tabular}{rrrrrrrrr}
& \multicolumn{2}{c}{Explicit} && \multicolumn{2}{c}{Semiexplicit} && \multicolumn{2}{c}{Implicit} \\
\multicolumn{1}{l}
{step size} & 
\multicolumn{1}{l}
{error} & 
\multicolumn{1}{l}
{order} & &
\multicolumn{1}{l}{error\rule{0pt}{2.5ex}\rule{0pt}{2.5ex}} 
& 
\multicolumn{1}{l}
{order} & &
\multicolumn{1}{l}{
error
\rule{0pt}{2.5ex}}
&
 \multicolumn{1}{l}
{order}\\
\hline

5.000e-02  & 4.864e-07  & --          & & 1.095e-08  & --     & & 4.606e-09  & --          \\

2.500e-02  & 3.250e-08  & 3.90       & & 2.952e-09  & 1.89   & & 2.919e-09  & 0.66       \\

1.250e-02  & 2.192e-09  & 3.89       & & 3.451e-10  & 3.10   & & 3.445e-10  & 3.08        \\

7.692e-03  & 3.281e-10  & 3.91       & & 6.201e-11  & 3.54   & & 6.196e-11  & 3.53        \\

4.545e-03  & 4.137e-11  & 3.94       & & 8.748e-12  & 3.72   & & 8.745e-12  & 3.72        \\

2.632e-03  & 4.757e-12  & 3.96       & & 1.075e-12  & 3.84   & & 1.075e-12  & 3.83        \\

1.563e-03  & 5.991e-13  & 3.97       & & 1.410e-13  & 3.90   & & 1.409e-13  & 3.90        \\
\label{exp1SDIRK}
\end{tabular}
\end{table}

\begin{table}[htpb]
\centering
\caption{Errors and orders of convergence of Example 1, $\lambda = 1$, Radau IA3}
\label{Exp7Results}
\begin{tabular}{rrrrrrrrr}
& 
\multicolumn{2}{c}
{Explicit} && 
\multicolumn{2}{c}
{Semiexplicit} && 
\multicolumn{2}{c}
{Implicit} \\
\multicolumn{1}{l}
{ step size} & 
\multicolumn{1}{l}
{error} & 
\multicolumn{1}{l}
{order} & &
\multicolumn{1}{l}
{error\rule{0pt}{2.5ex}\rule{0pt}{2.5ex}} & 
\multicolumn{1}{l}
{order} & &
\multicolumn{1}{l}
{error\rule{0pt}{2.5ex}} & 
\multicolumn{1}{l}
{order}\\
\hline

1.000e-01  & 6.995e-07  & --  & & 4.621e-08  & --  & & 3.550e-08  & --  \\

3.333e-02   & 2.858e-09  & 5.01  & & 1.587e-10  & 5.16  & & 1.553e-10  & 4.94  \\

2.000e-02  & 2.209e-10  & 5.01  & & 1.222e-11  & 5.02  & & 1.218e-11  & 4.98  \\

1.429e-02  & 4.094e-11  & 5.01  & & 2.274e-12  & 5.00  & & 2.274e-12  & 4.99  \\

1.111e-02   & 1.163e-11  & 5.01  & & 6.528e-13  & 4.97  & & 6.540e-13  & 4.96  \\

9.091e-03   & 4.261e-12  & 5.00  & & 2.378e-13  & 5.03  & & 2.377e-13  & 5.04  \\

\end{tabular}
\end{table}

As time integrators, we use the scheme (\ref{nonlinearscheme}) with the rational functions of the Runge--Kutta methods SDIRK3 ($p=4$) and the 3-stages RadauIA3 ($p=5$) \cite{Hairer}. Notice that the application of these RK methods does not give its classical order of convergence $p$. According to the main result in \cite{IC}, we should expect orders $p' = 3.25$ and $4.25$ for SDIRK3 and RadauIA3, while Tables \ref{exp1SDIRK} and \ref{Exp7Results} show that the scheme (\ref{nonlinearscheme}) avoids the order reduction, as it is predicted by Theorem~\ref{mainth}. Both tables show results taking $\lambda = 1$. Notice that both methods satisfy hypothesis 3 because they are $A$-stable and satisfy $|r_{\infty}|<1$.

In this case ($\lambda = 1$), it is interesting to note that in both cases the semiexplicit mode involves an improvement of the error by slightly more than an order of magnitude. However, the implicit mode does not practically improve on the semiexplicit mode, so its higher computational cost is not justified. Moreover, if we take the number of systems being solved in the integration as a magnitude for the computational cost, Figure \ref{fig1} suggests that the explicit mode is more efficient than the semiexplicit one, since the computational cost of the latter is twice that of the first for the same number of steps.

\begin{figure}[htpb]
\centering
\includegraphics[width=\linewidth]{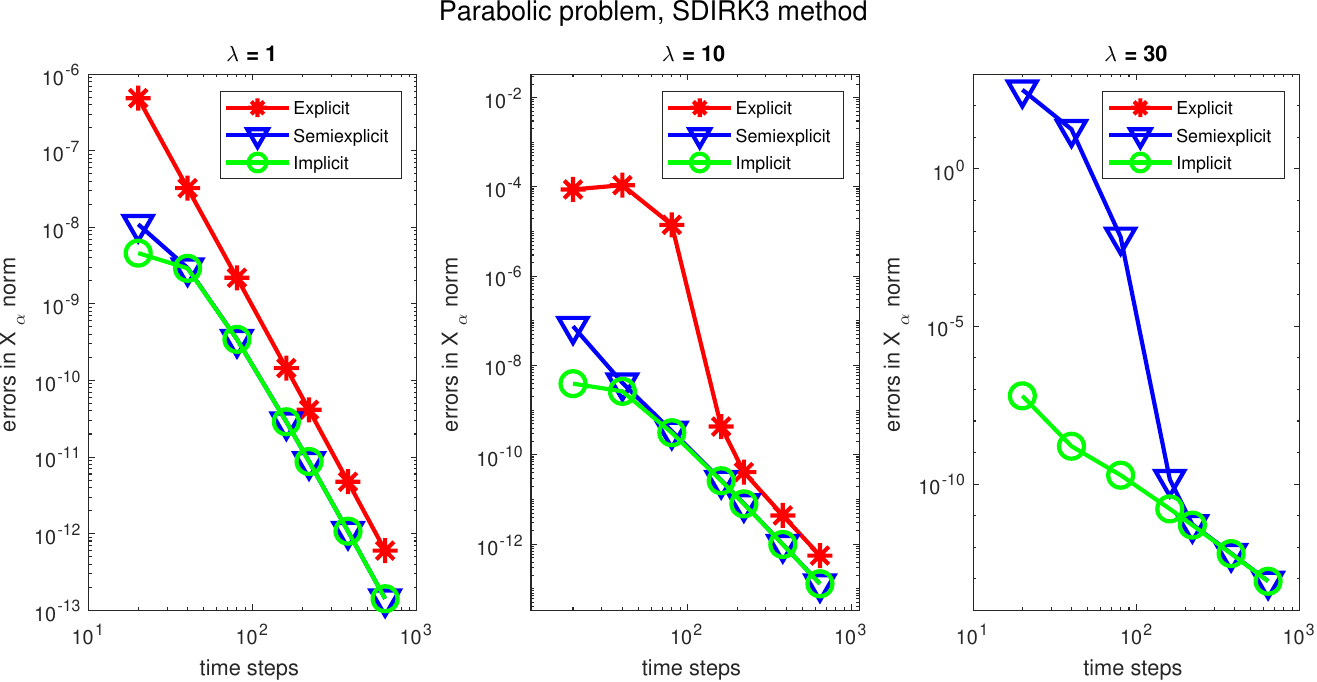}
\caption{This figure shows the error in the discrete norm $\|\cdot\|_{H^1(0,1)}$ for the parabolic problem with $\lambda = 1, 10, 30$ and the method SDIRK3.}
\label{fig2}
\end{figure}
\begin{figure}[htpb]
\centering
\includegraphics[width=\linewidth]{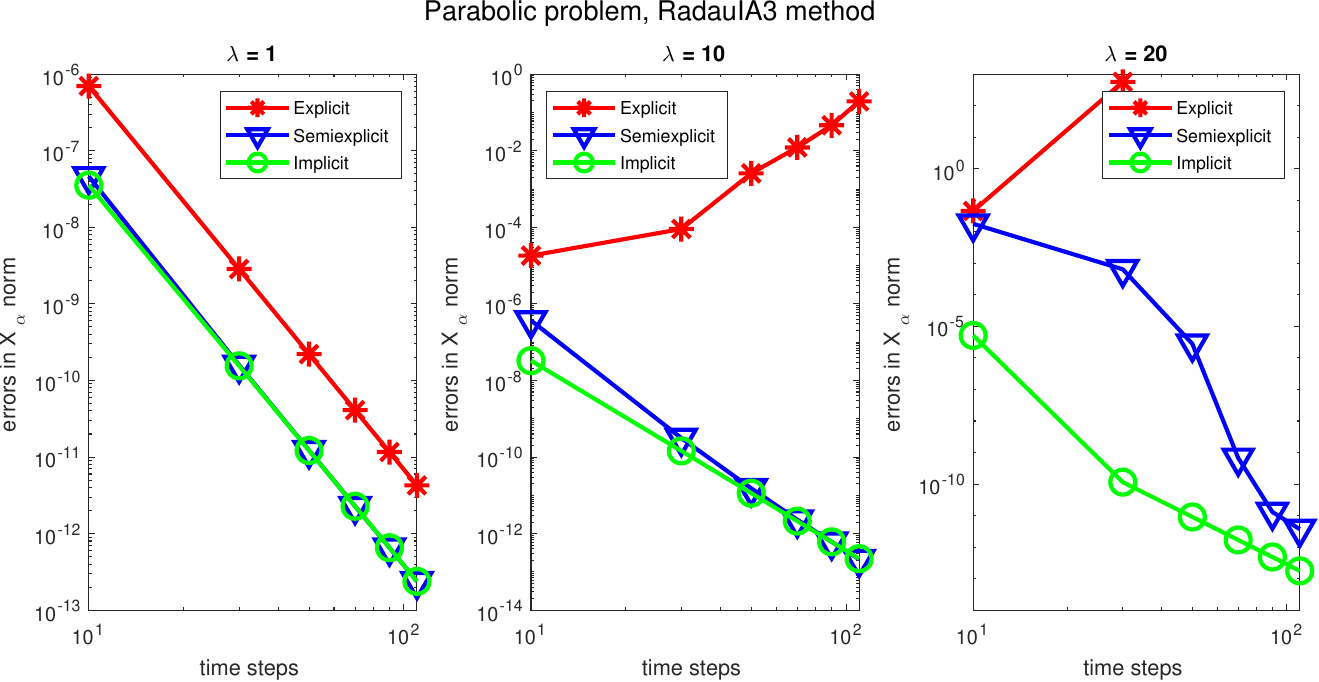}
\caption{This figure shows the error in the discrete norm $\|\cdot\|_{H^{3/2}(\Omega)}$ for the parabolic problem with $\lambda = 1, 10, 15$ and the method RadauIA3.}
\label{fig3}
\end{figure}
Figures \ref{fig2} and \ref{fig3} show the error in the integration when the parameter $\lambda$ is modified. We take $\lambda = 1, 10, 30$ for SDIRK3 and $\lambda = 1, 10, 15$ for RadauIA3. With both methods we observe that the explicit method does not perform out well when $\lambda$ increases, while the implicit one does not almost vary. The semiexplicit has an intermediate behaviour. This is consistent with the well-known sensitivity of explicit methods to stiff problems, whereas implicit methods handle it better due to their greater stability.

\paragraph*{Example 2.}  We consider a semilinear parabolic problem in the domain $\Omega = (0,1) \times (0,1)$ with homogeneous Dirichlet boundary conditions,
\begin{equation}
\left\{ \begin{array}{lcll} u_t(t,x,y) & = & \Delta u(t,x,y) + u^2 + f(t,x,y), & \quad 0\le t \le 1, \quad (x,y) \in \Omega, \\
u(0,x,y)  & = & u_0(x,y), & \quad (x,y) \in \Omega, \\
u(t,x,y) & = & 0, & \quad 0 \leq t \leq 1, \quad (x,y) \in 	\partial \Omega.
\end{array} \right.
\label{ej2}
\end{equation}
where $f : [0,1] \times \Omega \rightarrow \mathbb{C}$, $u_0 : \Omega \rightarrow \mathbb{C}$. In order to fit the problem in our framework, we take $X = L^2\left(\Omega\right)$, $A = \Delta$ with $D(A) = H^2\left(\Omega\right) \cap H^1_0\left(\Omega\right)$ and $\alpha = 3/4$, so that $\|\cdot\|_{\alpha}= \|\cdot\|_{H^{3/2}(\Omega)}$. We adjust the data $u_0$ and $f$ in such a way that $u(t,x,y) = x(1-x)y(1-y)\,\text{e}^t$, $0 \leq t, x \leq 1$, is the solution of the problem.

We discretize the problem in space by means of finite differences. To this end, we fix a number $J$ of uniformly distributed nodes $x_j = jh, 1 < j < J$, $y_k = kh, 1 < k < J$, in $(0, 1)$, with $h = 1/(J + 1)$. The spatial derivatives $u_{xx}$ and $u_{yy}$ are approximated by using the standard three-point finite difference scheme, respectively. In this way, since the exact solution is a polynomial of degree two in $x$ and $y$, there are no spatial errors.

\begin{table}[htpb]
\centering
\caption{Errors and orders of convergence of Example 2 with SDIRK3, $J=50$.}
\label{Exp10ResultsJ=50}
\begin{tabular}{rrrrrrrrr}
& 
\multicolumn{2}{c}
{Explicit} && 
\multicolumn{2}{c}
{Semiexplicit} && 
\multicolumn{2}{c}
{Implicit} \\
\multicolumn{1}{l}
{ step size} & 
\multicolumn{1}{l}
{error} & 
\multicolumn{1}{l}
{order} & &
\multicolumn{1}{l}
{error\rule{0pt}{2.5ex}\rule{0pt}{2.5ex}} & 
\multicolumn{1}{l}
{order} & &
\multicolumn{1}{l}
{error\rule{0pt}{2.5ex}} & 
\multicolumn{1}{l}
{order}\\
\hline

2.500e-02  & 1.402e-08  & --  & & 2.047e-10  & --  & & 1.791e-10  & --  \\

1.250e-02  & 9.234e-10  & 3.92  & & 9.012e-11  & 1.18  & & 8.930e-11  & 1.00  \\

6.250e-03  & 6.197e-11  & 3.90  & & 1.033e-11  & 3.12  & & 1.031e-11  & 3.11  \\

3.125e-03  & 4.095e-12  & 3.92  & & 8.581e-13  & 3.59  & & 8.571e-13  & 3.59  \\

1.563e-03  & 2.641e-13  & 3.95  & & 6.128e-14  & 3.81  & & 6.128e-14  & 3.81  \\

\end{tabular}
\end{table}

\begin{table}[htpb]
\centering
\caption{Errors and orders of convergence of Example 2 with RadauIA3, $J=50$.}
\label{Exp11Results}
\begin{tabular}{rrrrrrrrr}
& 
\multicolumn{2}{c}
{Explicit} && 
\multicolumn{2}{c}
{Semiexplicit} && 
\multicolumn{2}{c}
{Implicit} \\
{ step size} & 
\multicolumn{1}{l}
{error} & 
\multicolumn{1}{l}
{order} & &
\multicolumn{1}{l}
{error\rule{0pt}{2.5ex}\rule{0pt}{2.5ex}} & 
\multicolumn{1}{l}
{order} & &
\multicolumn{1}{l}
{error\rule{0pt}{2.5ex}} & 
\multicolumn{1}{l}
{order}\\
\hline

1.000e-01  & 5.600e-05  & --  & & 2.869e-06  & --  & & 3.188e-06  & --  \\

5.000e-02  & 2.322e-06  & 4.59  & & 1.241e-07  & 4.53  & & 1.320e-07  & 4.59  \\

2.500e-02  & 8.220e-08  & 4.82  & & 4.511e-09  & 4.78  & & 4.665e-09  & 4.82  \\

1.250e-02  & 2.721e-09  & 4.92  & & 1.515e-10  & 4.90  & & 1.542e-10  & 4.92  \\

6.250e-03  & 8.746e-11  & 4.96  & & 5.300e-12  & 4.84  & & 5.340e-12  & 4.85  \\

\end{tabular}
\end{table}

As time integrators, we use again the scheme (\ref{nonlinearscheme}) with the rational functions of the Runge--Kutta methods SDIRK3 ($p=4$) and 3-stages RadauIA3 ($p=5$). Notice that the application of these RK methods does not give its classical order of convergence $p$, while Tables \ref{Exp10ResultsJ=50} and \ref{Exp11Results} show that the scheme (\ref{nonlinearscheme}) does, as it is predicted by Theorem \ref{mainth}.

\paragraph*{Example 3.}  We consider a semilinear hyperbolic problem in the unit interval with periodic boundary conditions,
\begin{equation}
\left\{ \begin{array}{lcll} u_t(t,x) & = & - u_{x}(t,x) + u - u^3 + f(t,x), & \quad 0\le t \le 1, \quad 0 \leq x \leq 1, \\
u(0,x)  & = & u_0(x), & \quad 0 \leq x \leq 1, \\
u(t,0) & = & u(t,1), & \quad 0 \leq t \leq 1,
\end{array} \right.
\label{ej3}
\end{equation}
where $f : [0,1] \times [0,1] \rightarrow X$, $u_0 : [0,1] \rightarrow X$. In order to fit the problem in our framework, we take $X = H^1[0,1]$, $A = - d/dx$ with $D(A) = \left\lbrace u \in H^1[0,1] : u(0) = u(1) \right\rbrace$ and $\alpha = 0$, so that $\|\cdot\|_{\alpha}= \|\cdot\|_{X}=\|\cdot\|_{H^1(0,1)}$. We adjust the data $u_0$ and $f$ in such a way that $u(t,x) = x^3 \,\text{e}^t\,\sin(\pi\,x) + (1 -\text{e}^t)$, $0 \leq t, x \leq 1$, is the solution of the problem. We fix again number $J$ of uniformly distributed nodes $x_j = jh, 1 < j < J$, in $(0, 1)$, with $h = 1/(J + 1)$. The spatial derivatives are approximated by upwind finite differences.

Tables \ref{exp3SDIRK} and \ref{Exp8Results} show the results for the hyperbolic problem with the scheme (\ref{nonlinearscheme}) and the rational function of SDIRK3 and 3-stages RadauIA3, respectively. According to \cite{IC}, these RK methods have orders $p =3.5$ and $p=4.5$, respectively, when applied to this problem. We find the results to be similar to the parabolic one.

Again, we obtain an improvement on the size of the error after the semiexplicit correction for a fixed stepsize and we observe that the order of the methods is in agreement with that being predicted by Theorem \ref{mainth}.

\begin{figure}[htpb]
\centering
\includegraphics[height=70mm]{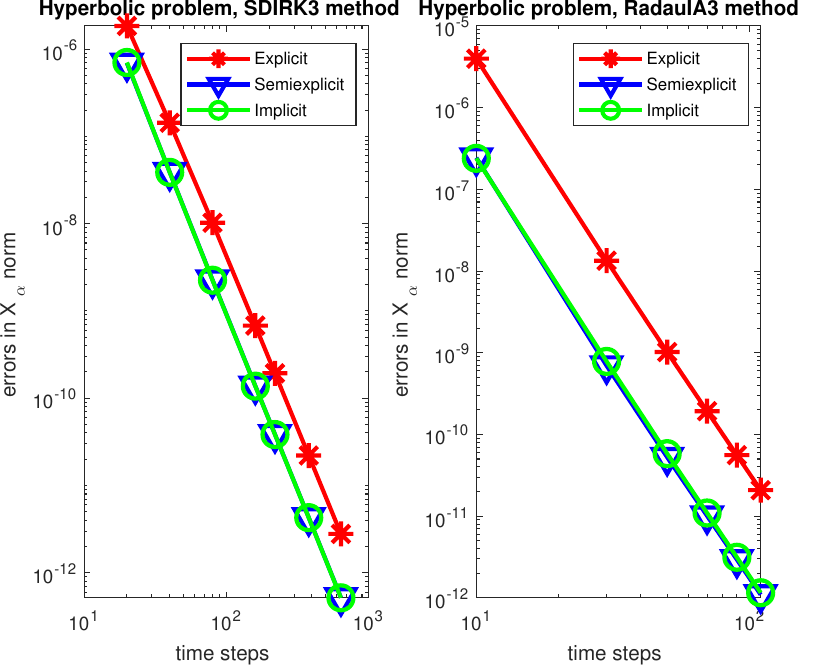}
\caption{This figure shows the error in the discrete norm associated to $\|\cdot\|_{H^1(0,1)}$ for the hyperbolic problem (\ref{ej3}) with two different time integrators and $J = 100$.}
\end{figure}

\begin{table}[ht]
\centering
\caption{Errors and order of convergence in Example 3 with the rational SDIRK3 method, J=100.}
\begin{tabular}{rrrrrrrrr}
& 
\multicolumn{2}{c}
{Explicit} && 
\multicolumn{2}{c}
{Semiexplicit} && 
\multicolumn{2}{c}
{Implicit} \\
\multicolumn{1}{l}
{ step size} & 
\multicolumn{1}{l}
{error} & 
\multicolumn{1}{l}
{order} & &
\multicolumn{1}{l}
{error\rule{0pt}{2.5ex}\rule{0pt}{2.5ex}} & 
\multicolumn{1}{l}
{order} & &
\multicolumn{1}{l}
{error\rule{0pt}{2.5ex}} & 
\multicolumn{1}{l}
{order}\\
\hline

5.000e-02  & 1.839e-06  & --          & & 6.896e-07  & --     & & 7.013e-07  & --          \\

2.500e-02  & 1.430e-07  & 3.69       & & 3.879e-08  & 4.15   & & 3.867e-08  & 4.18       \\

1.250e-02  & 1.024e-08  & 3.80       & & 2.266e-09  & 4.09   & & 2.251e-09  & 4.10        \\

7.692e-03  & 1.537e-19  & 3.91       & & 3.185e-10  & 4.04   & & 3.168e-10  & 4.04        \\

4.545e-03  & 1.929e-10  & 3.95       & & 3.849e-11  & 4.02   & & 3.837e-11  & 4.01        \\

2.632e-03  & 2.206e-11  & 3.97       & & 4.299e-12  & 4.01   & & 4.286e-12  & 4.01        \\

1.563e-03  & 2.791e-12  & 3.97       & & 5.206e-13  & 4.01   & & 5.200e-13  & 4.02        \\
\label{exp3SDIRK}
\end{tabular}
\end{table}

\begin{table}[htpb]
\centering
\caption{Errors and orders of convergence of Example 3 with RadauIA3, J=100.}
\label{Exp8Results}
\begin{tabular}{rrrrrrrrr}
& \multicolumn{2}{c}
{Explicit} && 
\multicolumn{2}{c}
{Semiexplicit} && 
\multicolumn{2}{c}
{Implicit} \\
\multicolumn{1}{l}
{ step size} & 
\multicolumn{1}{l}
{error} & 
\multicolumn{1}{l}
{order} & &
\multicolumn{1}{l}
{error\rule{0pt}{2.5ex}\rule{0pt}{2.5ex}} & 
\multicolumn{1}{l}
{order} & &
\multicolumn{1}{l}
{error\rule{0pt}{2.5ex}} & 
\multicolumn{1}{l}
{order}\\
\hline

1.000e-01  & 3.976e-06  & --  & & 2.440e-07  & --  & & 2.397e-07  & --  \\

3.333e-02  & 1.339e-08  & 5.18  & & 6.939e-10  & 5.34  & & 7.826e-10  & 5.21  \\

2.000e-02  & 1.018e-09  & 5.04  & & 5.255e-11  & 5.05  & & 5.815e-11  & 5.09  \\

1.429e-02  & 1.925e-10  & 4.95  & & 1.009e-11  & 4.91  & & 1.092e-11  & 4.97  \\

1.111e-02  & 5.582e-11  & 4.93  & & 2.942e-12  & 4.90  & & 3.135e-12  & 4.97  \\

9.091e-03  & 2.078e-11  & 4.92  & & 1.093e-12  & 4.93  & & 1.154e-12  & 4.98  \\

\end{tabular}
\end{table}

\section*{Funding}
The authors Carlos Arranz-Simón and Begoña Cano are supported by the Spanish Ministerio de Ciencia e Innovaci\'on under Project PID2023-147073NB-I00.



\end{document}